\documentclass[12pt]{amsart}
\usepackage{amssymb, amscd, amsmath, amsthm, latexsym, enumerate}
\newtheorem{theorem}{Theorem}
\newtheorem{lemma}[theorem]{Lemma}
\newtheorem*{cor}{Corollary}

\begin{document}

\title{Surgery on $\widetilde{\mathbb{SL}}\times\mathbb{E}^n$-manifolds}

\author[J. A. Hillman]{J. A. Hillman}
\address{School of Mathematics and Statistics\\ University of Sydney\\  
{Sydney}, NSW 2006\\ AUSTRALIA}
\email{jonh@maths.usyd.edu.au}

\author[S.K.Roushon]{S.K.Roushon}
\address{ School of Mathematics, Tata Institute, Mumbai 400-005, INDIA}
\email{roushon@math.tifr.res.in}

\begin{abstract}
We show that closed $\widetilde{\mathbb{SL}}\times\mathbb{E}^n$-manifolds 
(with $n\geq2$) are topologically rigid, if $n\geq2$, and are rigid up to
$s$-cobordism, if $n=1$.
\end{abstract}

\thanks{To appear in the Canadian Mathematical Bulletin}
\subjclass{57R67, 57N16}
\keywords{Topological rigidity, geometric structure,
surgery groups}

\maketitle


The ``Borel Conjecture" asserts that aspherical closed manifolds are
topologically rigid, i.e. are determined up to homeomorphism by their
fundamental groups.
Farrell and Jones have proven this for infrasolvmanifolds (of dimensions 
$\geq4$) \cite{FJ97} 
and for smooth manifolds (of dimensions $\geq5$)
with Riemannian metrics of nonpositive curvature \cite{FJ93},
while Nicas and Stark have shown it to hold for manifolds (of dimensions 
$\geq5$) admitting an effective codimension-2 toral action of hyperbolic type 
\cite{NS85}.
The work of Farrell and Jones was used in \cite{Ro00} to establish topological 
rigidity for $M\times D^k$ for all orientable closed irreducible 3-manifolds 
$M$ with $\beta_1(M)>0$ and all $k\geq3$.
In our main result (Theorem 4) we shall adapt the approach of \cite{Ro00} 
to show that manifolds with finite covering spaces admitting such toral actions 
are also topologically rigid.
(These are the natural higher dimensional analogues of 
geometric Seifert fibred spaces.)
In particular, all closed
$\widetilde{\mathbb{SL}}\times\mathbb{E}^n$-manifolds with $n\geq2$
are topologically rigid,
although such manifolds do not admit metrics of nonpositive curvature 
\cite{Eb82}, 
and may not admit codimension-2 toral actions.
(See Corollary A.) 

The most immediate applications of our work are to low dimensions.
In Corollary B we complete the characterization of aspherical geometric 
4-manifolds up to $s$-cobordism, in terms of fundamental group and Euler
characteristic.
(This was established for all other 4-dimensional geometries of aspherical type
in Theorem 9.11 of \cite{Hi}.)
In passing we also show (in Theorem 2)
that noncompact complete $\widetilde{\mathbb{SL}}$-manifolds of finite volume 
are homeomorphic to complete $\mathbb{H}^2\times\mathbb{E}^1$-manifolds 
of finite volume.
(This was stated without proof in \cite{Th}.)

If $G$ is a group let $G'$, $\zeta{G}$ and $\sqrt{G}$ denote the
commutator subgroup, centre and Hirsch-Plotkin radical of $G$,
respectively.
(In the cases considered below $\sqrt{G}$ is always abelian,
and thus is the unique maximal abelian normal subgroup.)
If $H$ is a subgroup of $G$ let $C_G(H)$ denote the centralizer of $H$
in $G$. (Thus $\zeta{G}=C_G(G)$.)
Let $E(n)=Isom(\mathbb{E}^n)=\mathbb{R}^n\rtimes{O(n)}$.

The following lemma is based on Lemma 9.5 of \cite{Hi}.
We give it here for the convenience of the reader.

\begin{lemma}
Let $\pi$ be a finitely generated group with normal subgroups $A\leq N$
such that $A$ is free abelian of rank $r$, $[\pi:N]<\infty$ and $N\cong
A\times{N/A}$. Then there is a homomorphism $f:\pi\to{E(r)}$ with image
a discrete cocompact subgroup and such that $f|_A$ is injective.
\end{lemma}

\begin{proof} 
Let $G=\pi/N$ and $M=N^{ab}\cong{A}\oplus(N/AN')$.
Then $M$ is a finitely generated $\mathbb{Z}[G]$-module and the
image of $A$ in $M$ is a $\mathbb{Z}[G]$-submodule. 
Extending coefficients to the rationals $\mathbb{Q}$ gives a natural inclusion 
$\mathbb{Q}A\leq\mathbb{Q}M$, since $A$ is a direct summand of $M$ 
(as an abelian group), 
and $\mathbb{Q}A$ is a $\mathbb{Q}[G]$-submodule of $\mathbb{Q}M$. 
Since $G$ is finite $\mathbb{Q}[G]$ is semisimple,
and so $\mathbb{Q}A$ is a $\mathbb{Q}[G]$-direct summand of $\mathbb{Q}M$.
Let $K$ be the kernel of the homomorphism from $M$ to $\mathbb{Q}A$
determined by a splitting homomorphism from $\mathbb{Q}M$ to $\mathbb{Q}A$,
and let $\tilde{K}$ be the preimage of $K$ in $\pi$.
Then $K$ is a $\mathbb{Z}[G]$-submodule of $M$ and $M/K\cong Z^r$,
since it is finitely generated and torsion free of rank $r$.
Moreover $\tilde K$ is a normal subgroup of $\pi$ and $A\cap\tilde{K}=1$.
Hence $H=\pi/\tilde{K}$ is an extension of $G$ by $M/K$ and
$A$ maps injectively onto a subgroup of finite index in $H$.
Let $T$ be the maximal finite normal subgroup of $H$. 
Then $H/T$ is isomorphic to a discrete cocompact subgroup of $E(r)$,
and the projection of $\pi$ onto $H/T$ is clearly injective on $A$.
\end{proof}

\begin{theorem}
Let $M$ be a $3$-manifold which is Seifert fibred over a complete open 
$\mathbb{H}^2$-orbifold $B$ of finite area.
Then $M$ is homeomorphic to a complete open 
$\mathbb{H}^2\times\mathbb{E}^1$-manifold.
\end{theorem}

\begin{proof}
Let $\pi=\pi_1(M)$ and let $A\cong Z$ be the image in $\pi$
of the fundamental group of the general fibre. 
Let $p:\pi\to\pi/A\cong\pi_1^{orb}(B)$ be the epimorphism 
given by the Seifert fibration,
and let $\psi:\pi_1^{orb}(B)\to Isom(\mathbb{H}^2)$ 
be a monomorphism onto a discrete subgroup of finite coarea
which determines the hyperbolic structure of $B$.

Since $B$ is complete and has finite area $\pi_1^{orb}(B)$ is finitely 
generated and since $B$ is open $\pi_1^{orb}(B)$ has a free normal subgroup 
$F$ of finite index.
Then $\pi$ is finitely generated.
Let $N=p^{-1}(F)\cap C_\pi(A)$.
Then $A<N$ and $N\cong A\times(N/A)$, 
since $A$ is central in $N$ and $N/A$ is free.
Hence there is a homomorphism $f:\pi\to E(1)$ which is injective on $A$, 
by the lemma.
Let $\theta=(\psi{p},f):\pi\to Isom(\mathbb{H}^2\times\mathbb{E}^1)$.
Then $\theta$ is injective, and $\theta(\pi)$ is a discrete subgroup 
of finite covolume.
Since $\theta(\pi)$ is torsion free it acts freely and so 
$N=H^2\times{R}/\theta(\pi)$ is a complete open 
$\mathbb{H}^2\times\mathbb{E}^1$-manifold of finite volume.
Projection from $H^2\times {R}$ onto the first factor induces 
a Seifert fibration of $N$ over $B$,
and since $\pi_1(N)\cong\pi=\pi_1(M)$ it follows that $M$ and $N$ are
homeomorphic.
\end{proof}

In particular, if $M$ is a compact 3-manifold with a nontrivial JSJ
decomposition then every geometric piece of type $\widetilde{\mathbb{SL}}$
also admits the geometry $\mathbb{H}^2\times\mathbb{E}^1$.
(This is part of Theorem 4.7.10 of \cite{Th}.
However no proof is given there.)

A similar argument shows that if $M$ is an open $(m+2)$-manifold which is the
total space of an orbifold bundle with base a complete open hyperbolic
2-orbifold $B$ of finite area, general fibre a flat $m$-manifold $F$ and 
monodromy group a finite subgroup of $Out(\pi_1(F))$ then there is an 
$\mathbb{H}^2\times\mathbb{E}^m$-manifold $M_1$ which is an orbifold bundle 
with base $B$ and general fibre $F$ and a homotopy equivalence $f:M\to M_1$ 
which preserves the conjugacy classes of the subgroups corresponding to 
the cusps.
Since the cusps are flat $(m+1)$-manifolds we may assume that $f$ is a
homeomorphism off a compact set, and a relative version of the Farrell-Jones
curvature argument then shows that $f$ is homotopic to a homeomorphism,
if $m\geq3$.
Is there a direct, elementary argument to show that $M$ and $M_1$ 
must be fibrewise diffeomorphic (for any $m\geq1$)?

If an $(m+2)$-manifold $M$ admits an effective action of the $m$-torus
$T^m=R^m/Z^m$ the image in $\pi_1(M)$ of the fundamental group of the 
principal orbit is central and the orbit space $Q$ is a 2-orbifold \cite{CR77}.
In \cite{NS85} it is shown that if $M$ is an aspherical $(m+2)$-manifold 
with an effective $T^m$ action of hyperbolic type the higher Whitehead 
groups $Wh_i(\pi_1(M))$ are trivial for all $i\geq0$ and  
$|S_{TOP}(M\times D^k,\partial)|=1$, whenever $m+k\geq4$
(or $m+k\geq3$, if $\partial M=\emptyset$).
Their argument for the Whitehead groups extends immediately to the following
situation.

\begin{lemma}
Let $\pi$ be a torsion-free group with a virtually poly-$Z$ normal subgroup $N$
such that $\pi/N\cong \pi_1^{orb}(B)$, where $B$ is a compact $2$-orbifold.
Then $Wh(\pi)=0$.
\end{lemma}

\begin{proof}
If $B$ is a closed $\mathbb{E}^2$-orbifold then $\pi$ is virtually poly-$Z$
and the result is proven in \cite{FH81}.
If $B$ is a closed $\mathbb{H}^2$-orbifold the argument of \cite{NS85} using
hyperelementary induction applies with little change.
If $\pi/N$ is virtually free it is the fundamental group of a 
graph of groups with all vertex groups finite or 2-ended and all edge groups
finite, and so $\pi$ is the fundamental group of a 
graph of groups with all vertex groups torsion free and virtually poly-$Z$.
Thus the result follows from \cite{FH81} and the Waldhausen Mayer-Vietoris
sequence
\cite{Wa78}.
(Note that $c.d.\pi<\infty$ since $\pi$ is torsion free, $c.d.N<\infty$ and 
$v.c.d.\pi/N\leq2$ in all cases.)
\end{proof}

The argument of \cite{NS85} determining the surgery structure sets for
$(m+2)$-manifolds admitting an effective $T^m$ action of hyperbolic type
appears to use the hypothesis of a toral action in an essential way, 
to establish an induction on $m$.
We shall rely instead on the curvature argument of \cite{FJ93}.

\begin{theorem}
Let $M$ be the total space of an orbifold bundle $p:M\to Q$ with base
$Q$ a closed $\mathbb{H}^2$-orbifold and general fibre a flat $m$-manifold $F$
of dimension $\geq3$, and such that $A=\sqrt{\pi_1(F)}\cong Z^m$ 
is centralized by a subgroup of finite index in $\pi=\pi_1(M)$.
If $f:M_1\to M$ is a homotopy equivalence with $M_1$ a closed $m$-manifold
then $f$ is homotopic to a homeomorphism.
\end{theorem}

\begin{proof}
Suppose first that there is an epimorphism $q:\pi_1^{orb}(Q)\to Z$.
Let $\hat{Q}$ and $\hat M$ be the induced covering spaces and $\hat p:\hat 
M\to \hat Q$ be the corresponding fiber bundle projection.
Then $\hat Q$ is noncompact, and is the increasing union 
$\hat{Q}=\cup_{k\geq1}{Q_k}$ of compact suborbifolds with nontrivial boundary.
We may assume that for each $k\geq0$ the boundary of $Q_k$ does not contain 
any corner points,
$G_k=\pi_1^{orb}(Q_k)$ is not virtually abelian,
and $G_k$ maps injectively to $G=\pi_1^{orb}(\hat{Q})$. 
Let $DQ_k$ be the closed orbifold obtained by doubling $Q_k$ along its boundary.
Since $\pi_1^{orb}(DQ_k)$ is not virtually abelian
there is a monomorphism $\psi:\pi_1^{orb}(DQ_k)\to Isom(\mathbb{H}^2)$
with image a discrete, cocompact subgroup.
(See page 248 of \cite{ZVC}.)

Let $M_k=\hat p^{-1}(Q_k)$. Then $M_k$ is a compact bounded $m$-manifold
and $\hat p:M_k\to Q_k$ is an orbifold fibration with general fibre $F$.
Doubling $M_k$ gives a closed $m$-manifold $DM_k$
with an orbifold fibration over $DQ_k$,
and $\pi(k)=\pi_1(DM_k)$ is an extension of $\pi_1^{orb}(DQ_k)$ by
$\pi_1(F)$.
As $\pi_1^{orb}(Q_k)$ acts on $A$ through a finite
subgroup the centralizer of $A$ in $\pi(k)$ 
again has finite index.
Let $N$ be a characteristic subgroup of finite index in $\pi(k)$ 
which centralizes $A$ and such that $N/A$ is a $PD_2^+$-group,
and let $e\in H^2(N/A;A)$ be the cohomology class of the 
extension $0\to{A}\to{N}\to{N/A}\to1$.
The reflection which interchanges the copies of $M_k$ leaves the boundary
pointwise fixed, and projects to the corresponding reflection of $DQ_k$.
Thus it induces an automorphism of $N$ which is the identity on $A$ and
reverses the orientation of $N/A$.
It follows that $e=-e$ and so the extension splits: $N\cong A\times{N/A}$.
Therefore there is a homomorphism $f:\pi(k)\to E(m)$ 
which is injective on $A$, by Lemma 1.
The homomorphism 
$(\psi_kp|_{\pi(k)},f):\pi(k)\to Isom(\mathbb{H}^2\times\mathbb{E}^m)$
has finite kernel, and so is injective, since $\pi$ is torsion free.
The quotient $P_k=H^2\times{R^m}/\pi(k)$ is closed and nonpositively curved,
and is Seifert fibred over $DQ_k$.
Moreover $DM_k\simeq P_k$ since each is aspherical,
and so $M_k$ is a homotopy retract of $P_k$.

Now the structure set of $P_k$ is trivial,
by the Topological Rigidity theorem of Farrell and Jones \cite{FJ93}. 
Since $M_k$ is a homotopy retract of $P_k$, 
the structure set of $M_k$ is also trivial. 
Equivalently, the assembly maps 
$H_j(M_k;\mathbb{L}_o^w)\to L_j(\pi_1(M_k),w)$ are isomorphisms for $j$ large,
where $w=w_1(M)$.
(Note that no decorations are needed on the surgery obstruction groups 
as $Wh(\pi)=0$, by Lemma 3.) 
Since homology and $L$-theory commute with direct limits we conclude that 
$H_j(\hat{M};\mathbb{L}_o^w)\to L_j(\pi_1(\hat{M}),w)$ is an isomorphism 
for $j$ large.
Using the Wang sequence for homology, 
naturality of the assembly maps and Ranicki's algebraic version of Cappell's
Mayer-Vietoris sequence for square root closed HNN extensions
it follows that the same is true for $M$.
(See \cite{Ro00} for more details.)

If $\beta_1(\pi_1^{orb}(Q))=0$ we may use hyperelementary induction,
as in \cite{NS85}, to reduce to the case already treated.
\end{proof}

A similar curvature argument could be used to show that $Wh(\pi)=0$,
for $\pi=\pi_1(M)$ as in the theorem.

\begin{cor}
\!{\bf A.}
Let $M$ be a closed $\widetilde{\mathbb{SL}}\times\mathbb{E}^n$-manifold,
where $n\geq2$, and let $f:M_1\to M$ be a homotopy equivalence,
with $M_1$ a closed $(n+3)$-manifold.
Then $f$ is homotopic to a homeomorphism.
\end{cor}

\begin{proof}
The composite of projection from the model space $\widetilde{SL}\times {R^n}$ 
onto the first factor with the fibration of $\widetilde{SL}$ over 
$\mathbb{H}^2$
induces an orbifold bundle fibration $p:M\to Q$, 
with base $Q$ a closed $\mathbb{H}^2$-orbifold and general fibre $F$ a flat 
$(n+1)$-manifold.
In Theorem 9.3 of \cite{Hi} it is shown that when $n=1$ the fundamental
group 
of a closed $\widetilde{\mathbb{SL}}\times\mathbb{E}^n$-manifold has a subgroup 
of finite index which is a direct product, 
and the argument extends immediately to the general case.
It follows that $\sqrt{\pi_1(F)}\cong Z^{n+1}$ 
is centralized by a subgroup of finite index in $\pi$,
and so we may apply the theorem. 
\end{proof}

We may adapt this result to obtain a somewhat weaker result for the case $n=1$ 
by taking products with $S^1$.

\begin{cor}
\!{\bf B.}
Let $N$ be a closed $\widetilde{\mathbb{SL}}\times\mathbb{E}^1$-manifold,
and let $N_1$ be a closed 4-manifold with $\pi_1(N_1)\cong\pi=\pi_1(N)$ 
and $\chi(N_1)=\chi(N)$.
Then $N_1$ is $s$-cobordant to $N$.
\end{cor}

\begin{proof}
The manifold $N_1$ is aspherical, 
by Corollary 3.5.1 of \cite{Hi},
and so there is a homotopy equivalence $g:N_1\to N$.
Let $M=N\times S^1$, $M_1=N_1\times S^1$, and $f=g\times{id}_{S^1}$.
Then $M$ is a $\widetilde{\mathbb{SL}}\times\mathbb{E}^2$-manifold.
Hence $f$ is homotopic to a homeomorphism $h:M_1\cong M$, by the theorem.
Since $h\sim{g\times{id}_{S^1}}$ it lifts to a homeomorphism 
$N_1\times{R}\cong{N}\times{R}$.
The submanifold of $N\times{R}$ bounded by $N\times\{0\}$ 
and a disjoint copy of $N_1$ is an $h$-cobordism. 
It is in fact an $s$-cobordism,
since $Wh(\pi_1(N))=0$, by Lemma 3.
\end{proof}

This result complements Theorem 9.11 of \cite{Hi},
where a similar result is proven for all 4-manifolds admitting a
nonpositively curved geometry.

Is there a corresponding result for manifolds with a proper geometric
decomposition?
The argument for Theorem 3.3 of \cite{Le95} extends readily to show that
if 
$M$ is a $n$-manifold with a finite collection of disjoint flat
hypersurfaces $\mathcal{S}$ such that the components of $M-\cup\mathcal{S}$
all have complete finite volume geometries of type $\mathbb{H}^n$ or
$\mathbb{H}^{n-1}\times\mathbb{E}^1$, and if
there is at least one piece of type $\mathbb{H}^n$
then $M$ admits a Riemannian metric of nonpositive sectional curvature
(see \cite{Eb82}).
Such manifolds are topologically rigid if $n\geq5$, by \cite{FJ93},
and we again deduce rigidity up to $s$-cobordism when $n=4$,
as in the above corollary.

\end{document}